\documentclass[12pt]{article}

\usepackage{geometry}
\usepackage[dvips]{graphicx}
\usepackage{float}
\usepackage{textcomp}
\usepackage{amsthm}
\usepackage{amsmath}
\usepackage{amssymb}
\usepackage{mathrsfs}
\usepackage{newtxmath}
\usepackage{booktabs}
\usepackage{rotating}
\usepackage{bm}
\usepackage{bbm}
\usepackage{tikz}

\usepackage[hidelinks]{hyperref}
\usepackage[backend=bibtex,maxnames=10]{biblatex}
\usepackage{color}
\usepackage{multirow}
\usepackage{tabularx}
\usepackage[all]{xy}
\usepackage{algorithm}
\usepackage{scalefnt}
\usepackage[shortlabels]{enumitem}
\usepackage{epstopdf}
\ifpdf
\DeclareGraphicsExtensions{.eps,.pdf,.png,.jpg}
\else
\DeclareGraphicsExtensions{.eps}
\fi
\newtheorem{theorem}{Theorem}[section]

\newtheorem{definition}[theorem]{Definition}

\usepackage{{algcompatible}}

\graphicspath{{FIG/}}
\usepackage{subfigure}

\newcommand{\R}{\mathbb{R}}

\newcommand{\set}[1]{\{#1\}}

\newcommand{\eps}{\eps}

\definecolor{orange}{rgb}{1.00,0.41,0.16}

\begin{document}

\title{Exploring polynomial models in the Search Step of Direct Multisearch}
\author{Ana L. Cust\'odio\thanks{Center for Mathematics and Applications (NOVA Math), Campus de Caparica, 2829-516 Caparica, Portugal \texttt{(alcustodio@fct.unl.pt)}. Support for this author was provided by Fundação para a Ciência e a Tecnologia (Portuguese Foundation for Science and Technology) under the scope of the projects UID/00297/2025 and UID/PRR/00297/2025.}\and Marta Pozzi \thanks{Computer Science Department, Bocconi University, Milan, Italy \texttt{(marta.pozzi2@phd.unibocconi.it)}} \and Everton J. Silva\thanks{Center for Mathematics and Applications (NOVA Math), Campus de Caparica, 2829-516 Caparica, Portugal \texttt{(ejo.silva@alumni.fct.unl.pt)}. Support for this author was provided by Fundação para a Ciência e a Tecnologia (Portuguese Foundation for Science and Technology) under the scope of the projects UID/00297/2025 and UID/PRR/00297/2025.}}
\date{}
\maketitle
\begin{abstract}
Direct Multisearch (DMS) is a class of direct-search algorithms
designed for multiobjective derivative-free optimization. Its
framework consists of an optional search step and a poll step, the
latter ensuring the corresponding theoretical convergence
properties. Recently, a search strategy based on the minimization
of quadratic polynomial models, constructed from previously
evaluated points, was proposed to improve the numerical efficiency
of the method. While the construction of these surrogate models
has been investigated, considerably less attention has been
devoted to how they should be jointly minimized, with min-max
scalarization typically being adopted.

This work investigates how different model-minimization strategies
influence the performance of DMS. To this end, alternative
strategies for exploring the quadratic polynomial models within
the search step of DMS are proposed and numerically assessed,
including one based on the recently proposed Improved Front
Steepest Descent algorithm.
\\

%
%
\end{abstract}

\section{Introduction}
Multiobjective optimization is a challenging scientific domain,
relevant both from theoretical and practical perspectives, in
which several conflicting objective functions must be optimized
simultaneously. Problems of this type arise in a wide range of
applications
\cite{Afshari_Hare_Tesfamariam_2019,Brito_Judice_2022,Handl_Kell_Knowles_2007,Jin_Sendhoff_2008,Potrebko_Fiege_Biagioli_Poleszczuk_2017,Rosso_Ciancio_DellOlmo_Salata_2020}.

In this work, we consider the following multiobjective
derivative-free optimization problem:
\begin{equation}
\begin{array}{l}
\min\ F(x)=\left(f_1(x),\ldots,f_q(x)\right)^\top \ \text{s.t.}\
x\in\Omega
\end{array}
\label{MOO}
\end{equation}
where $q\geq 2$, $\Omega \subseteq \R^n$ is the feasible region,
and $F:\Omega\subseteq \R^n\to(\R\cup\{+\infty\})^q$. It is
assumed that the derivatives of each $f_i$, $i\in
I:=\set{1,\ldots,q}$ are neither available nor can be reliably
approximated. This setting is common in simulation-based
optimization, where function evaluations arise from complex and
computationally expensive simulations. In such situations,
objective functions are often nonsmooth, noisy, or unreliable,
motivating the use of derivative-free methods. Comprehensive
reviews of single-objective derivative-free optimization methods
can be found
in~\cite{Audet_Hare_2017,Conn_Scheinberg_Vicente_2009}.

Standard approaches for multiobjective derivative-free
optimization include scalarization
techniques~\cite{Audet_Savard_Zghal_2008,Audet_Savard_Zghal_2010},
which transform the original problem into a sequence of
single-objective subproblems that can be addressed using classical
nonlinear optimization techniques.
BiMADS~\cite{Audet_Savard_Zghal_2008}, designed for biobjective
optimization, generates a sequence of scalarized reformulations
solved by MADS~\cite{Audet_DennisJr_2006}, exploiting the fact
that Pareto points can be ordered in the two-objective setting. To
extend this idea to problems with more than two objectives,
MULTIMADS was introduced in~\cite{Audet_Savard_Zghal_2010},
combining MADS with the Normal Boundary Intersection (NBI)
framework~\cite{Das_DennisJr_1998} to approximate the Pareto
front. Although these methods benefit from the robustness of MADS
in derivative-free settings, BiMADS is restricted to biobjective
problems, and a practical implementation of MULTIMADS is not
publicly available.

In~\cite{Cocchi_Liuzzi_Papini_Sciandrone_2018}, an unidirectional
method was proposed that generalizes implicit filtering to
bound-constrained multiobjective derivative-free optimization. At
each iteration, the algorithm constructs a simplex-based Jacobian
to approximate a multiobjective steepest-descent direction, which
is then explored via a line search. Line search strategies were
also considered in~\cite{Liuzzi_Lucidi_Rinaldi_2016}, for general
nonlinear constrained multiobjective problems, when proposing
DFMO, a derivative-free algorithm based on an exact-penalty
approach combined with sufficient-decrease line searches.

Unlike scalarization-based approaches, the methods proposed
in~\cite{Bigeon_LeDigabel_Salomon_2021,Cocchi_Liuzzi_Papini_Sciandrone_2018,Custodio_Madeira_Vaz_Vicente_2011,Liuzzi_Lucidi_Rinaldi_2016}
rely directly on the concept of Pareto dominance and avoid
aggregation schemes. In~\cite{Custodio_Madeira_Vaz_Vicente_2011},
the authors introduce Direct Multisearch (DMS), a generalization
of directional direct search to multiobjective derivative-free
optimization. DMS was the first practical derivative-free method
with theoretical convergence guarantees capable of approximating
the entire Pareto front, and it has served as the foundation for
subsequent developments.

Building on the ideas
in~\cite{Audet_Savard_Zghal_2010,Audet_Savard_Zghal_2008,Custodio_Madeira_Vaz_Vicente_2011},
DMultiMADS was proposed in~\cite{Bigeon_LeDigabel_Salomon_2021}, a
MADS-based extension for bound-constrained multiobjective
optimization. Unlike BiMADS and MULTIMADS, DMultiMADS does not
reduce the problem to a sequence of scalarized single-objective
subproblems. Instead, it maintains and updates a list of
nondominated points at each iteration, in the spirit of DMS.
Although DMultiMADS shares conceptual similarities with DMS, it
differs in the selection of incumbent points at each iteration,
which impacts its practical performance.

As a directional direct-search method, each iteration of
DMS~\cite{Custodio_Madeira_Vaz_Vicente_2011} follows the
search-poll paradigm, where the search step is optional and
convergence is ensured by the poll step.
In~\cite{Bras_Custodio_2020}, the authors introduced BoostDMS, a
variant of DMS that incorporates quadratic polynomial
interpolation and regression models into the search step to
enhance performance.

BoostDMS reuses previously evaluated points to construct local
quadratic surrogate models for each objective function, and its
search step relies on the successive minimization of combinations
of these surrogates. This approach extends to the multiobjective
setting the strategy originally developed
in~\cite{Custodio_Rocha_Vicente_2010} for single-objective
directional direct search, while preserving the theoretical
convergence properties of DMS. Since the search step determines
trial points evaluated by the algorithm, the procedure adopted for
jointly minimizing the surrogate models may have a significant
impact on the efficiency of the overall method.

Although considerable attention has been devoted to the
construction of surrogate models in derivative-free optimization,
the question of how these models should be exploited once
available has received comparatively little attention In DMS, this
choice directly determines trial points to be evaluated by the
algorithm and may therefore have a substantial influence on its
numerical performance.

In this work, different approaches for jointly minimizing the
surrogate models built during the search step of BoostDMS are
proposed, thereby defining alternative search strategies. In
particular, the $\epsilon$-constraint scalarization, the Normal
Boundary Intersection (NBI) framework, and a strategy based on the
recently introduced Improved Front Steepest Descent (IFSD)
algorithm~\cite{Lapucci_Mansueto_2023} are detailed and
numerically assessed as alternatives for model minimization within
the search step of DMS.


The paper is organized as follows. Sections 2 and 3 review Direct
Multisearch and model computation, respectively. Alternative
approaches for exploiting the quadratic polynomial models are
detailed in Section 4. Section 5 reports the numerical experiments
supporting the proposed search-step alternatives. Finally, Section
6 concludes the work with final remarks and directions for future
research.

\section{Direct Multisearch (DMS)}
\label{sec:DMS} Direct Multisearch (DMS), introduced
in~\cite{Custodio_Madeira_Vaz_Vicente_2011}, extends directional
direct search to multiobjective optimization. It directly exploits
the concept of Pareto dominance to generate approximations to the
complete Pareto front throughout the optimization process.

DMS is an algorithmic class encompassing several variants,
depending, for example, on the chosen globalization strategy or
the definition of the search step. Since its original
introduction, several improved versions have been proposed, namely
MultiGLODS~\cite{Custodio_Madeira_2018}, a global variant of DMS,
BoostDMS~\cite{Bras_Custodio_2020}, where a search step based on
quadratic polynomial models was introduced, its parallel version
Parallel
BoostDMS~\cite{Tavares_Bras_Custodio_Duarte_Medeiros_2023}, and
DMS-FILTER-IR~\cite{Silva_Custodio_2024}, developed for
constrained optimization through a filter approach combined with
an inexact restoration step.

The DMS framework has been successfully applied in several
contexts, including financial portfolio
optimization~\cite{Brito_Judice_2022,Brito_Sebastiao_Godinho_2016}
and the optimization of three-dimensional road
alignments~\cite{Hirpa_Hare_Lucet_Pushak_Tesfamariam_2016}, and
has proven competitive even when compared with derivative-based
multiobjective optimization
algorithms~\cite{Andreani_Custodio_Raydan_2022}.

This section presents a simplified algorithmic description of
BoostDMS and summarizes its main theoretical properties, inherited
from DMS. The following definition formalizes the concept of
Pareto dominance through the strict partial order induced by the
cone $\mathbb{R}^q_+=\{z\in\mathbb{R}^q\mid z\geq 0\}$.

\begin{definition}[Pareto dominance]
\label{def:pareto_dominance} Let $x, x'\in
\Omega\subseteq\mathbb{R}^n$. We say that $x$ dominates $x'$ if
$$F(x)~\prec_F~F(x'), \text{ i.e., if }
F(x')-F(x)\in\mathbb{R}^q_+\setminus\{\mathbf{0}\}.$$
\end{definition}

Pareto dominance is used to characterize both global and local
optimality. The efficient set is the subset of points in $\Omega$
that are not dominated by any other point in $\Omega$. The image
of the efficient set through the objective function defines the
Pareto front, which represents the solution of the multiobjective
optimization problem. A point $\overline{x}\in\Omega$ is said to
be a local Pareto minimizer of Problem~\eqref{MOO} if there exists
a neighborhood $\mathcal{N}$ of $\overline{x}$ in which
$\overline{x}$ is nondominated, that is, if there is no
$x\in\mathcal{N}\cap \Omega$ such that $F(x)\prec_F
F(\overline{x})$.

BoostDMS, similarly to DMS, can adopt two distinct globalization
strategies: one based on a sufficient decrease condition and
another based on the use of integer lattices. For simplicity, and
since it corresponds to the available numerical implementation,
only the latter will be described here, meaning that all points
generated by the algorithm belong to an implicit mesh. For a more
general description, we refer the reader to the original
work~\cite{Custodio_Madeira_Vaz_Vicente_2011}.

Throughout the iterations, BoostDMS maintains and updates a list
of feasible nondominated points and corresponding stepsize
parameters, defining the current approximation to the Pareto
front. Constraints are handled through an extreme barrier
approach, meaning that only feasible points are evaluated in the
objective function, while infeasible points are assigned the value
$(+\infty,\ldots,+\infty)\in\mathbb{R}^q$. An iteration is
declared successful whenever at least one new feasible
nondominated point is found. Such points are added to the list,
which is then updated by removing all dominated points.

Each iteration of BoostDMS is divided into four steps. First, a
point and corresponding stepsize parameter are selected from the
list of feasible nondominated points as the current iterate. Any
selection criterion may be adopted, since it does not affect
convergence. In the original implementation of DMS, as well as in
BoostDMS, the selection relies on a spread metric, designed to
reduce gaps between consecutive points in the list.

An optional Search Step follows, constituting the distinguishing
feature of BoostDMS. Previously evaluated points are reused to
construct local quadratic polynomial models for each objective
function component. These points are selected within a ball
centered at the current iterate, with radius proportional to the
current stepsize parameter. The resulting surrogate models are
then minimized, generating new trial points for evaluation in the
original objective function.

Depending on the number of available sample points, the quadratic
models may be underdetermined (handled through minimum
Frobenius-norm interpolation), determined (corresponding to
polynomial interpolation), or overdetermined (constructed through
regression techniques). During the initial iterations, the search
step is typically skipped until a predefined sample threshold is
reached.

The models are minimized through a level-based procedure,
individually or in combination, inside a trust-region
corresponding to the ball used for sample selection. At the first
level, each individual quadratic model is minimized separately. If
among the resulting points there isn't a feasible nondominated
point of the original problem, the algorithm proceeds to the
second level, where pairs of surrogate models are combined through
a $\min$--$\max$ scalarization and minimized again within the
trust-region.The procedure terminates as soon as at least one
feasible nondominated point is found or all levels have been
exhausted. If the search step produces at least one new feasible
nondominated point, the iteration is declared successful and the
poll step is skipped. Otherwise, the poll step becomes mandatory.

In the Poll Step, a local search is performed around the current
iterate by evaluating the objective function at points generated
along a set of directions scaled by the current stepsize
parameter. The requirements imposed on the set of directions are
directly related to convergence conditions. In general, the
directions should locally conform to the geometry of the feasible
region and be asymptotically dense in the unit sphere, being
typically positive-spanning sets. For a detailed discussion, we
refer the reader to the original
work~\cite{Custodio_Madeira_Vaz_Vicente_2011}. If a new feasible
nondominated point is found during the poll step, the iteration is
declared successful. Otherwise, the list remains unchanged and the
iteration is deemed unsuccessful.

At the end of each iteration, the stepsize parameter is updated
for the current iterate (if it remains in the list) and for all
newly generated points. If the iteration is successful, the
stepsize parameter is maintained or increased. Otherwise, it is
decreased.

The overall structure of this simplified description of BoostDMS
is outlined in Algorithm~\ref{DMS}.

\begin{algorithm}[htp!]
\begin{description}
{\small \vspace{1ex}
\item[Initialization] \ \\
Choose $x_0\in \Omega$ such that $f_i(x_0)< +\infty$, $\forall
i\in \set{1,\ldots,q}$, $\alpha_0 > 0$ an initial stepsize
parameter, $0 < \beta_1\leq \beta_2 < 1$ the coefficients for
stepsize contraction and $\gamma \geq 1$ the coefficient for
stepsize expansion. Let $\mathcal{D}$ be a (possibly infinite) set
of positive spanning sets, with directions $d$ satisfying
$0<d_{\min}\leq\|d\|\leq d_{\max}$. Consider
$L_0=\set{(x_0;\alpha_0)}$ the initial list of nondominated points
and corresponding stepsize parameters.\vspace{1ex}

\item[For $k=0,1,2,\ldots$] \ \\
\begin{enumerate}
\item[1.] \textbf{Selection of an iterate point:} Order the list $L_k$ according to some criteria and select $(x_k;\alpha_k)\in L_k$ as the current iterate and stepsize parameter.\ \\

\item[2.] \textbf{Search step:} By reusing previously evaluated
points, build a quadratic polynomial model for each component of
the objective function, $f_i,\; i\in\{1,\ldots,q\}$. Minimize
combinations of these models, computing a finite set of points
$P_s=\{z_s\mid s\in S\}$. Evaluate $F$ at
$L_{\texttt{add}}=P_s\cap\Omega$. Let $L_{\texttt{trial}}$ be the
set obtained from $L_k\cup \{(x;\alpha_k):x \in
L_{\texttt{add}}\}$ by removing dominated points. If
$L_{\texttt{trial}} \neq L_k$, then declare the iteration (and the
search step) successful, set
$L_{k+1}=L_{\texttt{trial}}$, and go to Step 4.\ \\

\item[3.] \textbf{Poll step:} Choose a positive spanning set $D_k$
from the set $\mathcal{D}$ and consider the set
$P_k=\{x_k+\alpha_k d\mid d\in D_k\}$. Evaluate $F$ at
$L_{\texttt{add}}= P_k\cap \Omega$. Let $L_{\texttt{trial}}$ be
the set obtained from $L_k\cup \{(x;\alpha_k):x \in
L_{\texttt{add}}\}$ by removing dominated points. If
$L_{\texttt{trial}} \neq L_k$, then declare the iteration (and the
poll step) as successful and set $L_{k+1}=L_{\texttt{trial}}$.
Otherwise, declare the iteration as unsuccessful and set
$L_{k+1}=L_k$.\ \\

\item[4.] \textbf{Stepsize parameter update:} If the iteration was
successful, then maintain or increase the corresponding stepsize
parameter, by considering $\alpha_{k,\texttt{new}}\in
[\alpha_k,\gamma\alpha_k]$. Replace all the new points $(x_k
+\alpha_k d; \alpha_k)$ in $L_{k+1}$ by $(x_k+ \alpha_k d;
\alpha_{k,\texttt{new}})$, when success is coming from the poll
step, or $(z_s;\alpha_k)$ in $L_{k+1}$ by
$(z_s;\alpha_{k,\texttt{new}})$, when success is coming from the
search step. Replace also $(x_k;\alpha_k)$, if in $L_{k+1}$, by
$(x_k;\alpha_{k,\texttt{new}})$. Otherwise, decrease the stepsize
parameter, by choosing $\alpha_{k,\texttt{new}}\in
[\beta_1\alpha_k,\beta_2\alpha_k]$, and replace the poll pair
$(x_k;\alpha_k)$ in $L_{k+1}$ by $(x_k;\alpha_{k,\texttt{new}})$.
\vspace{1ex}
\end{enumerate}
\item[EndFor]\ \\}
\end{description}
\caption{A simplified description of BoostDMS} \label{DMS}
\end{algorithm}

Being an algorithmic variant of DMS, BoostDMS inherits all
associated convergence results. The convergence analysis begins by
establishing the existence of a subsequence of stepsizes
converging to zero. Convergence is then studied along particular
sequences of unsuccessful iterates, referred to as refining
subsequences.

\begin{definition}
\label{def:refining_subsequence} A subsequence $\{x_k\}_{k\in K}$
of iterates generated by Algorithm~\ref{DMS}, corresponding to
unsuccessful poll steps, is said to be a refining subsequence if
$\{\alpha_k\}_{k\in K}$ converges to zero. A limit point
$\overline{x}\in \R^n$ of a refining subsequence is called a
refined point.
\end{definition}

The behavior of the algorithm is subsequently analyzed through
limit directions associated with refined points, called refining
directions.

\begin{definition}
\label{def:refining_directions} Let $\bar{x}$ be the refined point
of a convergent refining subsequence $\{x_k\}_{k\in K}$. If
$\lim_{k\in K'} \frac{d_k}{\|d_k\|}$ exists, where $K'\subseteq
K$, $d_k\in D_k$, and $x_k+\alpha_k d_k \in \Omega$, for
sufficiently large $k \in K'$, then this limit is said to be a
refining direction for~$\bar{x}$.
\end{definition}

Assuming that the set of refining directions is asymptotically
dense in the Clarke tangent cone to $\Omega$ at refined points,
convergence to Pareto-Clarke critical points or Pareto-Clarke-KKT
critical points can be established.

\begin{theorem}\cite[Thm. 4.9]{Custodio_Madeira_Vaz_Vicente_2011}
Consider a refining subsequence $\{x_k\}_{k\in K}$ converging to
$\bar{x}\in\Omega$. Assume that $F$ is Lipschitz continuous near
$\bar{x}$ and $int(T_{\Omega}^{Cl}(\bar{x}))\neq\emptyset$. If the
set of refining directions for $\bar{x}$ is dense in
$T_\Omega^{Cl}(\bar{x})$, then $\bar{x}$ is a Pareto-Clarke
critical point, i.e, $$\forall d\in T_\Omega^{Cl}(\bar{x}),\exists
\ell=\ell(d)\in\{1,\ldots,q\} \text{ such that }
f_{\ell}^\circ(\bar{x};d) \geq 0.$$ If, in addition, $F$ is
strictly differentiable at $\bar{x}$, then $\bar{x}$ is a
Pareto-Clarke-KKT critical point, i.e.,  $$\forall d\in
T_\Omega^{Cl}(\bar{x}),\exists \ell=\ell(d)\in\{1,\ldots,q\}
\text{ such that } \nabla f_{\ell}(\bar{x})^\top d \geq 0.$$
\end{theorem}

\section{Quadratic polynomial models}
\label{sec:quad_models} Since the work of
Powell~\cite{MJDPowell_1964,MJDPowell_1970}, quadratic polynomial
interpolation models have become an important tool in
derivative-free optimization, acting as local surrogates of the
objective function. Whether used within a trust-region framework
or in the definition of a search step for single-objective
derivative-free optimization, they naturally emerged in BoostDMS
as models for the different components of the objective function.
In this section, we summarize the main ideas and properties
associated with these models, adapting the presentation
in~\cite{Bras_Custodio_2020} to our notation and framework.

Let $f : \mathbb{R}^n \to \mathbb{R}$ be a smooth function for
which derivatives are unavailable. Then $f$ can be approximated in
a neighborhood of a current iterate $x_k$ by a polynomial model $m
: \mathbb{R}^n \to \mathbb{R}$ of the form
$$m(y) = \xi^{\top} \phi(y),$$ where $\phi(y)$ is a polynomial
basis of dimension $s+1$ and $\xi$ is a vector of coefficients.
Given a set of previously evaluated sample points $Y = \{ y_0,
y_1, \ldots, y_p \}$ and  corresponding function values $f(Y) =
(f(y_0), \ldots, f(y_p))^{\top}$, with $y_0=x_k$, the model
coefficients are obtained by solving the linear system $ M(\phi,
Y)\, \xi = f(Y)$, where
$$
M(\phi, Y) = \left [ \begin{array}{cccc}\phi_0(y_0) & \phi_1(y_0) & \cdots & \phi_{s}(y_0) \\
\vdots & \vdots & \vdots & \vdots\\
\phi_0(y_{p}) & \phi_1(y_{p}) & \cdots & \phi_{s}(y_{p})
\end{array}\right ]
\;\;\;
f(Y) = \left[ \begin{array}{c} f(y_0) \\
\\ \vdots \\ f(y_p) \end{array} \right]
$$
This system is square when the number of interpolation points
equals the dimension of the chosen polynomial basis. In the
remaining cases, namely when the system is underdetermined or
overdetermined, the model coefficients can be computed through a
minimum-norm solution or least-squares
regression~\cite{Bras_Custodio_2020}, respectively.

For simplicity, let us consider the natural basis of monomials.
Under standard smoothness assumptions, namely the Lipschitz
continuity of the gradient, even underdetermined models can
recover the accuracy of first-order Taylor expansions, provided
that the sample set satisfies a suitable geometric property, known
as $\Lambda$-poisedness~\cite{Conn_Scheinberg_Vicente_2009}.

\begin{theorem}
Let $f:\mathbb{R}^n\to\mathbb{R}$ be continuously differentiable
and assume that its gradient is Lipschitz continuous with constant
$C_{\nabla f}>0$ on the ball
$$B(x_k;\Delta_k)=\{\,x\in\mathbb{R}^n\mid \ \|x-x_k\|<\Delta_k\,\}.$$
Let $Y=\{y_0,\ldots,y_p\}$ be a sample set with $n+1\leq
p<\frac{(n+1)(n+2)}{2}-1$, which is $\Lambda_L$-poised for linear
interpolation (or linear regression) at $x_k$, and let $m$ be a
quadratic model built from $Y$. Then there exists a constant
$C_p>0$ (depending only on $p$) such that, for all $y\in
B(x_k;\Delta_k)$,
\begin{eqnarray*}
&& \|\nabla f(y)-\nabla m(y)\| \;\leq\; C_p\,\Lambda_L\left(C_{\nabla f}+\|H\|\right)\,\Delta_k,\\
&& \left|\,f(y)-m(y)\,\right| \;\leq\;
\left(C_p\,\Lambda_L+\tfrac{1}{2}\right)\left(C_{\nabla
f}+\|H\|\right)\,\Delta_k^{2},
\end{eqnarray*}
with $H$ denoting the Hessian of the model.
\end{theorem}

These relations show that the interpolation error decreases
proportionally to the trust-region radius, provided that the norm
of the model Hessian remains controlled, motivating the Minimum
Frobenius Norm (MFN) approach for computing minimum-norm
solutions.

Formally, let $g$ and $H$ denote the gradient and Hessian of the
model, respectively. The MFN model is computed as the solution of
\begin{eqnarray*}
\min &&  \tfrac{1}{4}\,\| H \|_F^2\\
\text{s.t.} && f(y_i) = f(y_0) + g^{\top}(y_i - y_0) +
\tfrac{1}{2}(y_i - y_0)^{\top} H (y_i - y_0) \quad i=1,\ldots,p,
\end{eqnarray*}
where $\|\cdot\|_F$ denotes the Frobenius
norm~\cite{Bras_Custodio_2020}. This construction balances the
number of required sample points and the ability of the model to
incorporate curvature information.

When the number of sample points is sufficient to build a complete
quadratic model, Taylor-like error bounds can also be established
for the Hessian approximation. Again, an appropriate geometric
condition must be satisfied by the sample
set~\cite{Conn_Scheinberg_Vicente_2009}.

\begin{theorem}
Let $f$ be a twice continuously differentiable function with a
Lipschitz continuous Hessian (with Lipschitz constant $C_{\nabla^2
f}
> 0$) in the ball $B(x_k; \Delta_k)$. If the sample set $Y =
\{y_0, \ldots, y_p\}$, with $p +1 \geq \frac{(n+1)(n+2)}{2}$, is
$\Lambda_Q$-poised for quadratic interpolation (or quadratic
regression), then for all $y\in B(x_k;\Delta_k)$
\begin{eqnarray*}
\|\nabla^2 f(y) - \nabla^2 m(y)\| &\leq& C_p^1 \Lambda_Q C_{\nabla^2 f} \Delta_k,\\
\|\nabla f(y) - \nabla m(y)\| &\leq& C_p^2 \Lambda_Q C_{\nabla^2 f} \Delta_k^2,\\
| f(y) - m(y) | &\leq& \left( C_p^3 \Lambda_Q C_{\nabla^2 f} +
\frac{C_{\nabla^2 f}}{6} \right) \Delta_k^3,
\end{eqnarray*}
where $C_p^i$, $i \in \{1, 2, 3\}$, are positive constants
depending on $p$.
\end{theorem}

In~\cite{Custodio_Rocha_Vicente_2010}, the use of quadratic
polynomial models, built from previously evaluated points, was
proposed for the definition of a search step in single-objective
direct search. BoostDMS extends this strategy to direct
multisearch, for multiobjective optimization. Previously evaluated
points are again selected within a ball centered at the current
iterate and used to construct a model for each objective function
component. The type of model considered (underdetermined,
determined, or overdetermined) depends on the number of available
points. In all cases, the radius of the ball used for sample
selection is directly related to the stepsize parameter, ensuring
the quality of the models associated with refining subsequences,
since the corresponding sequence of stepsize parameters converges
to zero. For more details, the original references could be
consulted~\cite{Conn_Scheinberg_Vicente_2009,Custodio_Rocha_Vicente_2010,Bras_Custodio_2020}.

The remaining question is how to exploit these models in the
definition of a search step for DMS. Solving the corresponding
multiobjective minimization problem, namely

\begin{equation}
\begin{array}{l}
\min\ M_k(x)=\left(m^k_1(x),\ldots,m^k_q(x)\right)^\top \
\text{s.t.}\ x\in\Omega, \|x-x_k\|\leq \Delta_k
\end{array}
\end{equation}

\noindent requires a careful approach. Although this is a
derivative-based problem, the solution strategy should be selected
judiciously. Should an approximation to the Pareto front be
sought? Should scalarization approaches be considered? If so,
which ones?

In~\cite{Bras_Custodio_2020}, a weighted Chebyshev norm
scalarization is applied to the joint minimization of models for a
subset of components of the objective function by solving the
problems:

\begin{equation} \label{eq:MOOsearch}
\begin{array}{rl}
\min & \zeta\\[1ex]
{\rm s.t.} & m^k_i(x) \; \le \; \zeta, \quad i \in I\\
& \|x-x_k \| \; \le \; \Delta_k \\
&  x \in \Omega,
\end{array}
\end{equation}
where $I \subseteq\{1,2,\ldots,q\}$.

The remainder of the paper investigates how the choice of the
model-minimization strategy influences the performance of
BoostDMS, thereby assessing whether this procedure should itself
be regarded as an algorithmic design component rather than as a
secondary implementation detail.

\section{Alternative approaches for using quadratic models}
In this section, three different strategies for exploiting the
quadratic polynomial models within the search step of BoostDMS are
proposed. Specifically, we consider the Normal Boundary
Intersection (NBI) method~\cite{Das_DennisJr_1998}, an adapted
$\epsilon$-constraint approach, and the Improved Front Steepest
Descent (IFSD) algorithm~\cite{Lapucci_Mansueto_2023}. These three
approaches were selected because they represent fundamentally
different philosophies for solving multiobjective optimization
problems. The adapted $\epsilon$-constraint method follows a
scalarization approach, NBI relies on a geometric construction
designed to improve the distribution of solutions, whereas IFSD
directly exploits multiobjective descent directions. This
diversity allows us to investigate whether the manner in which the
surrogate models are minimized has a systematic influence on the
performance of BoostDMS.

\subsection{The Normal Boundary Intersection (NBI) method}
The Normal Boundary Intersection (NBI) method, originally
introduced in~\cite{Das_DennisJr_1998}, is a scalarization
approach designed to generate a well-distributed approximation of
the Pareto front in multiobjective optimization. NBI constructs a
family of subproblems that project candidate solutions along
directions normal to the so-called \emph{convex hull of individual
minima} (CHIM) in the objective space. This geometric construction
mitigates the uneven spacing of points typically associated with
classical scalarization approaches and is also capable of
revealing nonconvex regions of the Pareto front.

Let $M(x) = (m_1(x), m_2(x), \ldots, m_q(x))^\top$ denote the
vector of objective functions, and let $x_i^* = \arg\min_{x \in
\Omega} m_i(x)$, for $i = 1, \ldots, q$, be the individual global
minimizers of each objective. The corresponding vector $M^* =
(m_1(x_1^*), \ldots, m_q(x_q^*))^\top$, often referred to as the
\emph{ideal point}, serves as basis from which the payoff matrix
$\Phi$ is defined as $$ \Phi = \left[ M(x_1^*) - M^*, \; M(x_2^*)
- M^*, \; \ldots, \; M(x_q^*) - M^* \right].$$ The CHIM is then
described by the set  $\{\, M^* + \Phi \beta \mid \beta \in \Delta
\,\}$, where $\Delta = \{\beta \in \mathbb{R}^q_+\mid \sum_i
\beta_i = 1 \}$ is the unit simplex.  For each chosen weight
vector $\beta$, a unit normal vector $n$ to the CHIM is computed,
and the classical NBI subproblem is given by
\begin{eqnarray}
\label{eq:nbi-original}
\max_{x \in \Omega,\, t \in \mathbb{R}} && t \\
\text{s.t.} && M(x) - M^* = \Phi \beta + t n.\nonumber
\end{eqnarray}
The optimal solution of~\eqref{eq:nbi-original} corresponds to the
intersection between the feasible image of $M(x)$ and the ray
departing from the CHIM in direction $n$, yielding a
Pareto-optimal solution for the chosen $\beta$.

The formulation above may generate spurious Pareto-front
approximations, particularly when the search region is nonconvex.
Indeed, the method aims at generating boundary points, rather than
necessarily Pareto-optimal ones. To address this issue, a modified
version of the problem was proposed in~\cite{PKShukla_2007}:
\begin{eqnarray}
\label{eq:nbi-model}
\max_{x \in \Omega,\, t \in \mathbb{R}} && t \\
\text{s.t.} && M(x) - M^* \leq \Phi \beta + t\, n.\nonumber
\end{eqnarray}
Although relaxed, this formulation preserves the geometric essence
of NBI by driving the search along the normal direction, while
improving robustness under model uncertainty.

In BoostDMS, the NBI mechanism is embedded within the model-based
search step that precedes the poll step. Each quadratic model is
minimized individually, over the feasible region, subject to the
trust-region constraint, to estimate the local minima $x_i^*$. The
ideal point $M^*$ and model-based payoff matrix $\Phi$ are then
computed from these approximate minimizers. The CHIM and its
associated normal direction $n$ are subsequently constructed from
these model-based quantities.

To generate additional trial points, $2^q-1$ equidistant weight
vectors on the simplex are considered, matching the maximum number
of trial points that can be generated during the search step in
the current implementation of BoostDMS (corresponding to the
different models combinations inside each level). The $2^q-1$
modified NBI subproblems~\eqref{eq:nbi-model} are solved for the
selected weight vectors respecting the original constraints
defining $\Omega$ and within the trust-region. This ensures that
the resulting trial points remain in a region where the surrogate
models provide reliable approximations.

\subsection{Adapted $\epsilon-$constraint method}
The $\epsilon$-constraint method is a classical scalarization
technique for multiobjective optimization. It transforms a
multiobjective problem into a sequence of constrained
single-objective problems, where one objective is minimized while
the remaining objectives are restricted by upper bounds, denoted
by $\epsilon_j$. By systematically varying these bounds, the
method can generate representative Pareto-optimal solutions, even
when the Pareto front is nonconvex.

Given $M(x) = (m_1(x), \ldots, m_q(x))^\top$ and the feasible
region $\Omega$, the standard $\epsilon$-constraint problem
associated with the $i$-th objective is formulated as
\begin{eqnarray}
\label{eq:eps-standard}
\min_{x \in \Omega} && m_i(x) \\
\text{s.t.} && m_j(x) \leq \epsilon_j, \quad j = 1, \ldots, q, \;
j \neq i.\nonumber
\end{eqnarray}

Each choice of the bounds $\epsilon_j$ defines a distinct
subproblem, whose minimizer corresponds to a different trade-off
point on the Pareto front. This approach is particularly effective
when one objective is prioritized while the remaining ones act as
constraints controlling the allowed degradation.

In BoostDMS the classical formulation~\eqref{eq:eps-standard} is
adapted to operate on local quadratic surrogate models, within the
trust-region framework.

Let $x_k$ denote the current iterate, and let $m_i^k(x)$ represent
the quadratic polynomial model associated with the objective
function $f_i(x)$, constructed from sample points selected in the
ball centered at $x_k$ with radius $\Delta_k$. The
$\epsilon$-bounds are dynamically fixed at the current objective
values, namely $\epsilon_j = f_j(x_k)$ for $j \neq i$, so that the
search focuses on improving one objective without significantly
deteriorating the remaining ones. The resulting model-based
subproblem is given by
\begin{equation}
\label{eq:eps-model}
\begin{aligned}
\min_{x \in \Omega} \quad & m_i^k(x) \\
\text{s.t.} \quad & m_j^k(x) \leq \epsilon_j = f_j(x_k), \quad j \neq i, \\
& \|x_k-x\| \leq \Delta_k.
\end{aligned}
\end{equation}

This formulation ensures that the generated trial point remains
within a region where the surrogate models are expected to provide
accurate approximations of the objective function components. By
solving~\eqref{eq:eps-model} for each objective
$i\in\{1,\ldots,q\}$, a set of $q$ candidate points is obtained,
which are then evaluated using the true objective function $F$ in
an attempt to identify new feasible nondominated points.

\subsection{Improved front steepest descent (IFSD)}
The Improved Front Steepest Descent (IFSD)
method~\cite{Lapucci_Mansueto_2023} builds upon the classical
steepest descent algorithm for multiobjective
optimization~\cite{Fliege_Svaiter_2000}. However, in the spirit of
DMS, IFSD explicitly maintains a list of nondominated points that
is iteratively refined to approximate the Pareto front. Each
iteration combines \emph{refinement steps}, which drive the points
toward Pareto stationarity, with \emph{exploration steps}, which
enrich the front by identifying new trade-offs among the
objectives.

Given the smooth multiobjective problem
\begin{equation}
\label{mop_unc} \min_{x\in \mathbb{R}^n} M(x) = (m_1(x), \ldots,
m_q(x))^\top,
\end{equation}
the steepest-descent direction is defined by $$v(x) =
\arg\min_{d\in \mathbb{R}^n} \max_{j=1,\ldots,q} \left\{ \nabla
m_j(x)^\top d + \tfrac{1}{2}\|d\|^2 \right\}.$$ For any subset of
objectives $I \subset \{1,\ldots,q\}$, the corresponding
\emph{partial descent direction} is defined as $$v_I(x) =
\arg\min_{d \in \mathbb{R}^n} \max_{j \in I}\left\{ \nabla
m_j(x)^\top d + \tfrac{1}{2}\|d\|^2 \right\}.$$

\noindent The optimal values of these subproblems, denoted by
$\theta(x)$ and $\theta_I(x)$, respectively, measure the local
potential for common or partial descent. Pareto-stationary points
satisfy $\theta(x)=0$.

At each iteration of the IFSD
algorithm~\cite{Lapucci_Mansueto_2023}, every nondominated point
$x_c$ in the current list undergoes a two-stage procedure
consisting of a \emph{refinement step} followed by an
\emph{exploration step}.

In the refinement step, a common steepest descent direction
$v(x_c)$ is computed, and an Armijo-type line search is performed
to obtain
\[
z_c = x_c + \alpha_c v(x_c),
\]
ensuring sufficient decrease in all objectives, whenever
$\theta(x_c)<0$. This step guarantees convergence toward Pareto
stationarity.

The exploration step aims at densifying the approximation to the
Pareto front. Starting from $z_c$, for every subset $I \subset
\{1,\ldots,q\}$ such that $\theta_I(z_c) < 0$, a partial descent
direction $v_I(z_c)$ is computed. Line searches along these
directions generate candidates $z_c + \alpha_I v_I(z_c)$,
improving subsets of objectives. Newly generated points are added
to the list, while dominated solutions are removed. This mechanism
ensures that all points may contribute to the exploration process.

Under standard smoothness and compactness assumptions, the
resulting sequences of sets converge to sets of Pareto-stationary
points, and every accumulation point of the iterates is Pareto
stationary for Problem~\eqref{mop_unc}.


IFSD was originally designed for unconstrained multiobjective
optimization. However, the surrogate models computed during the
search step of BoostDMS are only considered reliable within the
trust-region $\|x-x_k\|\leq\Delta_k$. Therefore, the computation
of the descent direction and the partial descent directions took
into consideration the trust-region and the original problem
constraints. A different possibility, also tested but with worst
numerical performance, was to solve the initial unconstrained
problems and from the final approximation to the Pareto front of
Problem~\eqref{mop_unc}, only keep the points feasible with
respect to both $\Omega$ and the trust-region constraint.

Even so, this set is typically large. Considering the expensive
nature of the true objective function, only a small subset of
these points can be selected for evaluation. The model values are
assumed to provide accurate approximations of the objective
function components, and only points that remain nondominated with
respect to the current list $L_k$ (using model values as
surrogates for the true objective function values) are retained.
Let $L_k^{\tt IFSD}$ denote this reduced list of points.

If the list $L_k^{\tt IFSD}$ contains at most $2^q-1$ points, then
all of them are selected for evaluation, that is, $P_s=L_k^{\tt
IFSD}$. Otherwise, a maximum of $2^q-1$ points must be selected
from $L_k^{\tt IFSD}$, prioritizing the most isolated ones in an
attempt to densify the current approximation to the Pareto front.

\begin{figure}
\centering
\includegraphics[scale=0.8]{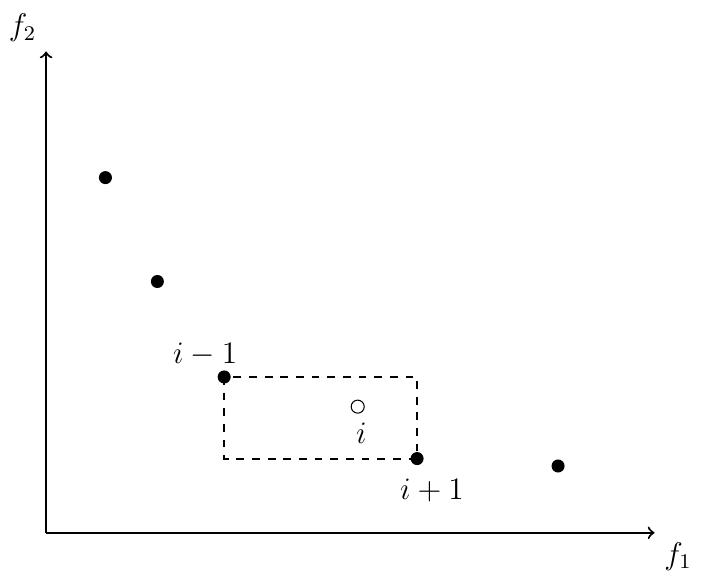}
\caption{Example of the computation of crowding distance and
centroid.} \label{crowd_centroide}
\end{figure}

To identify these isolated candidates, the crowding
distance~\cite{Deb_Pratap_Agarwal_Meyarivan_2002} is employed. For
each candidate point, the crowding distance is computed with
respect to an auxiliary list $L_{\tt aux}$, initially defined as
$L_k$, still using the model values as surrogates for the
objective function values at the new points. Points associated
with the least crowded regions of the objective space are
therefore preferred. In the case of ties, the selected point is
the one closest to the centroid of the corresponding crowding box
(see Figure~\ref{crowd_centroide}). Once selected, the point is
added to $P_s$ and included in $L_{\tt aux}$, again using the
model values as surrogates for the objective function values. The
selection strategy is repeatedly applied until a maximum of $2^q -
1$ points have been selected for evaluation in the true objective
function.

Figure~\ref{crowd_centroide} illustrates the concept of crowding
distance in the biobjective case. The black dots represent the
points in $L_{\tt aux}$ in the objective space $(f_1,f_2)$, while
the white dot corresponds to the model solution indexed by $i$ in
the final set computed by IFSD for model minimization. The
neighboring points, denoted by $i-1$ and $i+1$, define the bounds
of a rectangular region associated with point $i$, represented by
the dashed box. This region is constructed using the objective
values of adjacent points in each objective direction. The
crowding distance quantifies the extent of this region through
normalized differences between the objective values of neighboring
points. It therefore provides an estimate of the local density of
solutions around point $i$, where larger values indicate that the
solution lies in a less crowded region of the objective space,
thereby promoting diversity among the selected candidates. The
centroid of this box is defined as the midpoint of the rectangle
determined by the neighboring points $i-1$ and $i+1$.

\section{Numerical experiments}
\label{sec:numerical_exp} In this section, our goal is to evaluate
the numerical performance of the proposed search strategies.
Comparisons were made between the original BoostDMS and the three
proposed schemes. All tests were conducted on a server equipped
with 512 GB of RAM and two Intel{\small\textregistered}\
Xeon{\small\textregistered}\ Gold 6534 processors, each featuring
8 cores and 16 threads, with a base clock speed of 3.9 GHz, and a
maximum turbo frequency of 4.2 GHz. The algorithms were
implemented and executed using MATLAB R2024b.

\subsection{Performance assessment}
To compare the proposed strategies, we employ performance
profiles, introduced by Dolan and Mor\'e~\cite{Dolan_More_2002},
which provide a comprehensive framework for benchmarking
optimization algorithms across a collection of test problems. A
performance profile represents the cumulative distribution of a
given performance metric, allowing to simultaneously assess
efficiency and robustness. For a solver $s\in S$ and a set of
problems $P$, the profile is defined by
\begin{equation*}
\rho_s(\tau)=\frac{1}{|P|}\left|\left\{p\in P \mid t_{p,s}\leq \tau \min_{s\in S} t_{p,s} \right\}\right|,
\end{equation*}
where $\tau\geq1$ and $t_{p,s}$ denotes the value of a selected
performance metric obtained by solver $s$ on problem $p$. Thus,
$\rho_s(\tau)$ corresponds to the fraction of problems for which
solver $s$ performs within a factor $\tau$ of the best solver. In
particular, $\rho_s(1)$ measures the proportion of problems for
which solver $s$ is the most efficient, whereas larger values of
$\rho_s(\tau)$ for increasing $\tau$ indicate greater robustness.

The numerical comparison was carried out using four widely adopted
quality indicators: purity, hypervolume, and the spread measures
$\Gamma$ and $\Delta$. The purity metric quantifies the proportion
of nondominated points generated by a solver that belong to the
reference Pareto front. Let $F_{p,s}$ denote the approximation to
the Pareto front produced by solver $s$ for problem $p$, and let
$F_p$ be a reference Pareto front obtained by merging the
approximations generated by all solvers for the same problem and
removing dominated
points~\cite{Custodio_Madeira_Vaz_Vicente_2011}. The purity value
is then defined as
$$\bar{t}_{p,s}=\mathrm{Pur}_{p,s} =\frac{|F_{p,s}\cap F_p|}{|F_{p,s}|}.$$

To assess both convergence and diversity, we also considered the
hypervolume
metric~\cite{Zitzler_Thiele_Laumanns_Fonseca_Fonseca_2003}, which
measures the volume of the objective space dominated by the
approximation $F_{p,s}$ with respect to a reference point
$U_p\in\mathbb{R}^q$, dominated by all the points generated by the
different solvers for the same problem:
$$\bar{t}_{p,s}=\mathrm{HV}_{p,s}=\mathrm{Vol}\left( \bigcup_{x\in F_{p,s}}[x,U_p]\right),$$
where $\mathrm{Vol}(\cdot)$ denotes the $q$-dimensional Lebesgue
measure and $[x,U_p]$ is the hyper-rectangle defined by the lower
corner $x$ and the upper corner $U_p$.

Since larger values of purity and hypervolume indicate better performance, the corresponding performance profiles were constructed using their reciprocals, that is, $ t_{p,s}=\frac{1}{\bar{t}_{p,s}},$ so that lower values consistently represent superior performance across all metrics.

To complement the assessment of the Pareto front distribution, we
further considered the spread metrics $\Gamma$ and
$\Delta$~\cite{Deb_Pratap_Agarwal_Meyarivan_2002}. The $\Gamma$
metric measures the largest gap between consecutive nondominated
points, whereas the $\Delta$ metric evaluates the uniformity of
their distribution.

Suppose that solver $s$ computes for problem $p$ an approximation
consisting of the points $x^1,\ldots,x^N$. Let $x^0$ and $x^{N+1}$
denote the extreme points, corresponding to the points with the
best and worst values for each component of the objective
function~\cite{Custodio_Madeira_Vaz_Vicente_2011}. Then, the
$\Gamma$ metric is given by
\begin{equation}
\Gamma_{p,s}=\max_{j\in\{1,\ldots,q\}}\left(\max_{i\in{0,\ldots,N}}\delta_{j,i}\right),
\end{equation}
where $\delta_{j,i}=f_j(x^{i+1})-f_j(x^i),$ assuming that the objective values are sorted in ascending order for each objective component.

The $\Delta$ metric is defined as
\begin{equation}
\Delta_{p,s}=\max_{j\in\{1,\ldots,q\}} \left(
\frac{\delta_{j,0}+\delta_{j,N}+\sum_{i=1}^{N-1}|\delta_{j,i}-\bar{\delta}_j|}{\delta_{j,0}+\delta_{j,N}+(N-1)\bar{\delta}_j}
\right),
\end{equation}
where $\bar{\delta}_j$ denotes the average of the distances
$\delta_{j,i}$, $i=1,\ldots,N-1$. Lower values of both $\Gamma$
and $\Delta$ indicate better performance.

The assessment was conducted on a comprehensive benchmark
consisting of 100 bound-constrained test problems proposed
in~\cite{Custodio_Madeira_Vaz_Vicente_2011}, covering a broad
range of objective landscapes and difficulty levels. Regarding the
parameters of BoostDMS, all default values were considered. For
all variants, we set $\alpha_0=1$, $\beta_1=\beta_2=0.5$,
$\gamma=1$, and allowed a maximum of $5000$ function evaluations,
jointly with a minimum step size of $10^{-3}$, for all the points
in the list.

\subsection{The NBI method}
\begin{figure}[H]
\centering
\subfigure[Purity]{\includegraphics[scale=0.25]{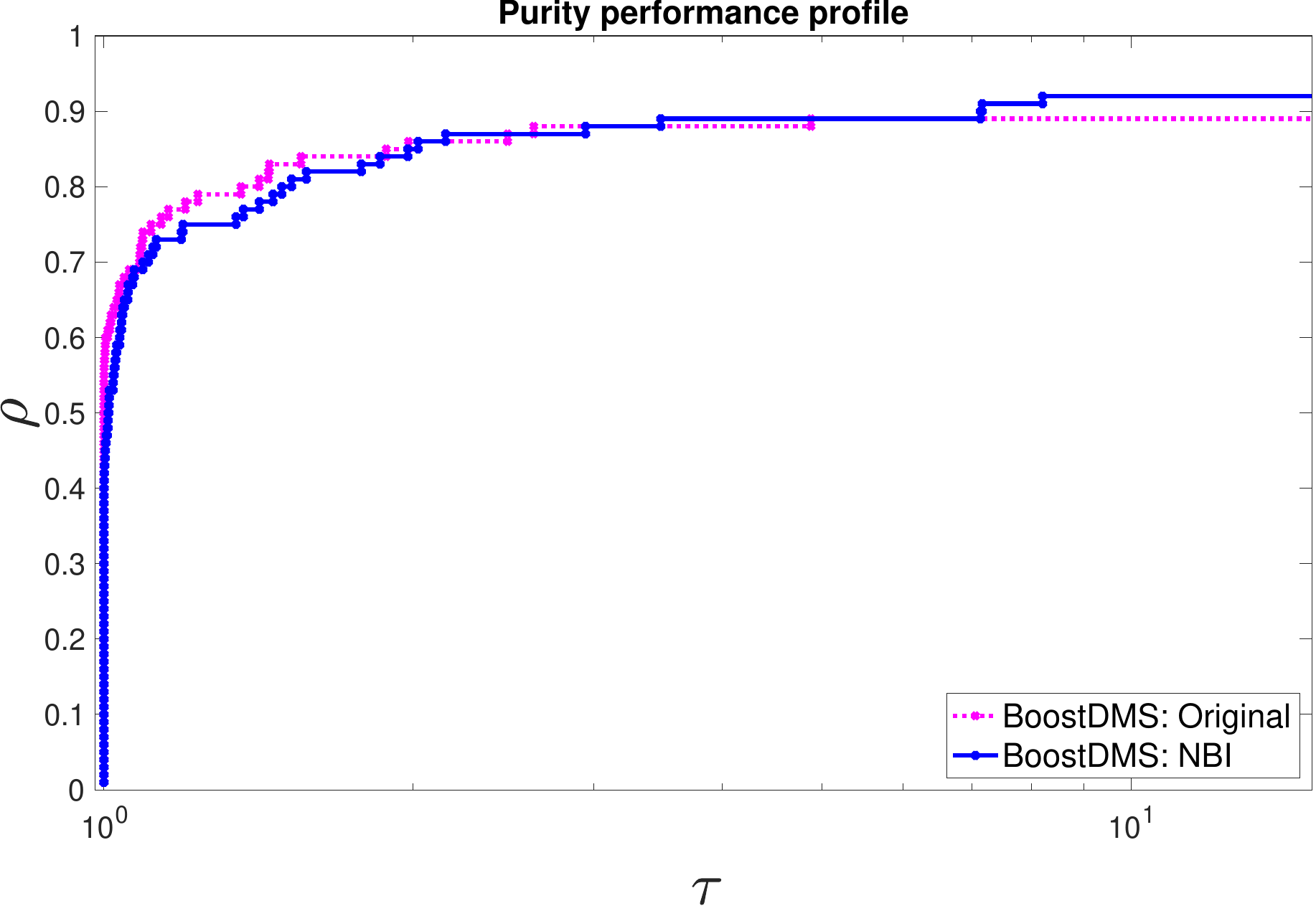}}
\subfigure[Hypervolume]{\includegraphics[scale=0.25]{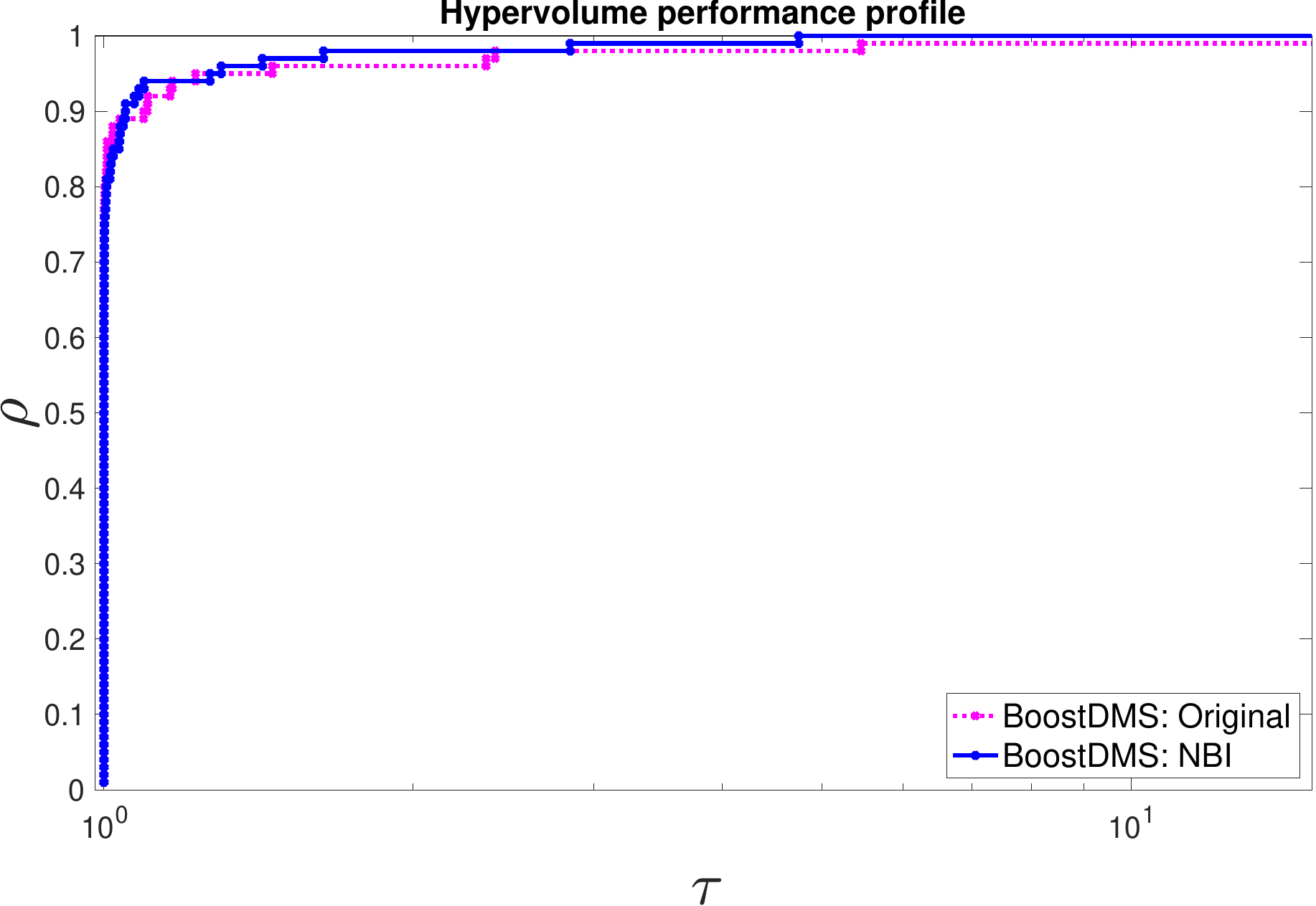}}
\subfigure[Spread Gamma
]{\includegraphics[scale=0.25]{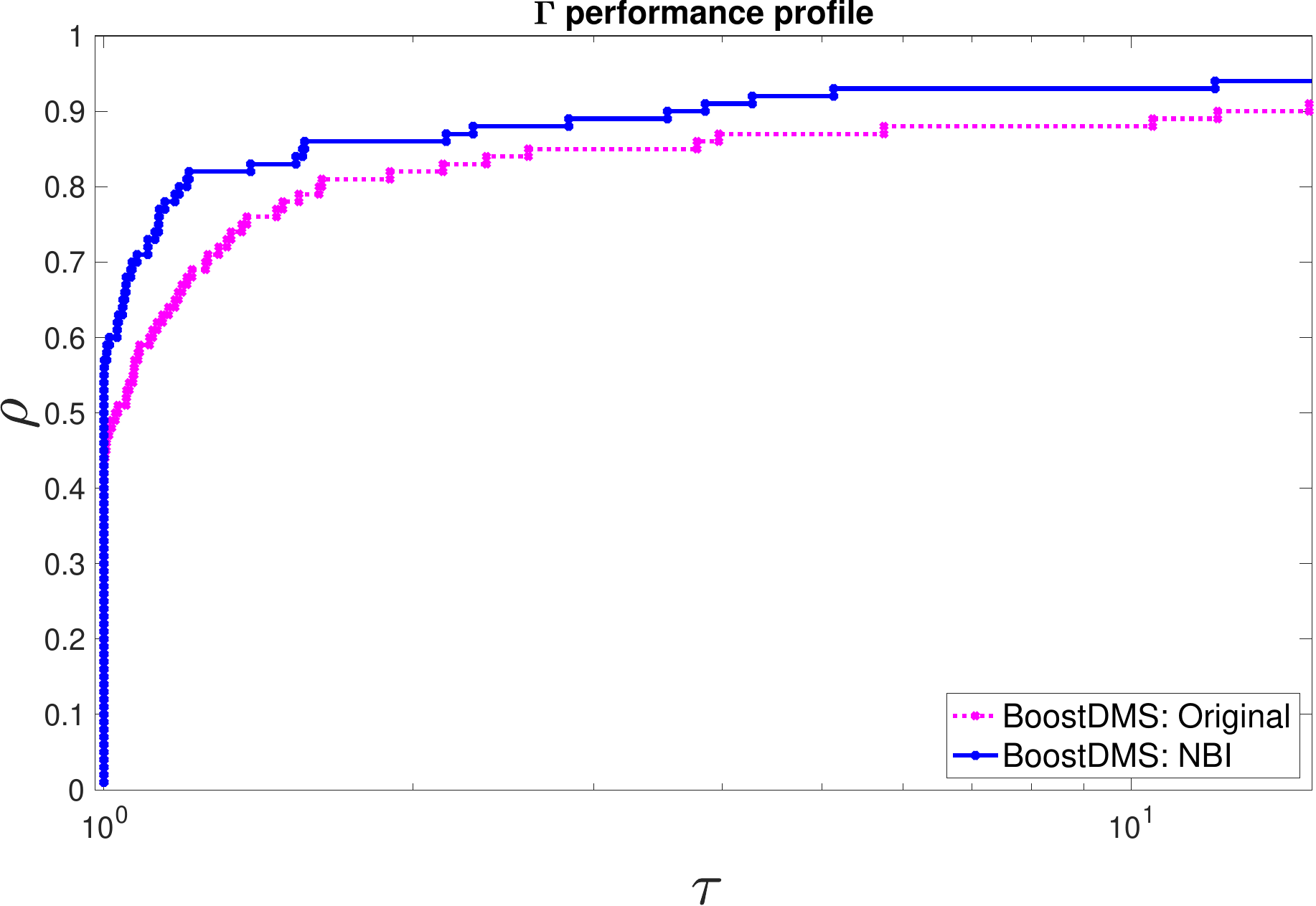}}
\subfigure[Spread
Delta]{\includegraphics[scale=0.25]{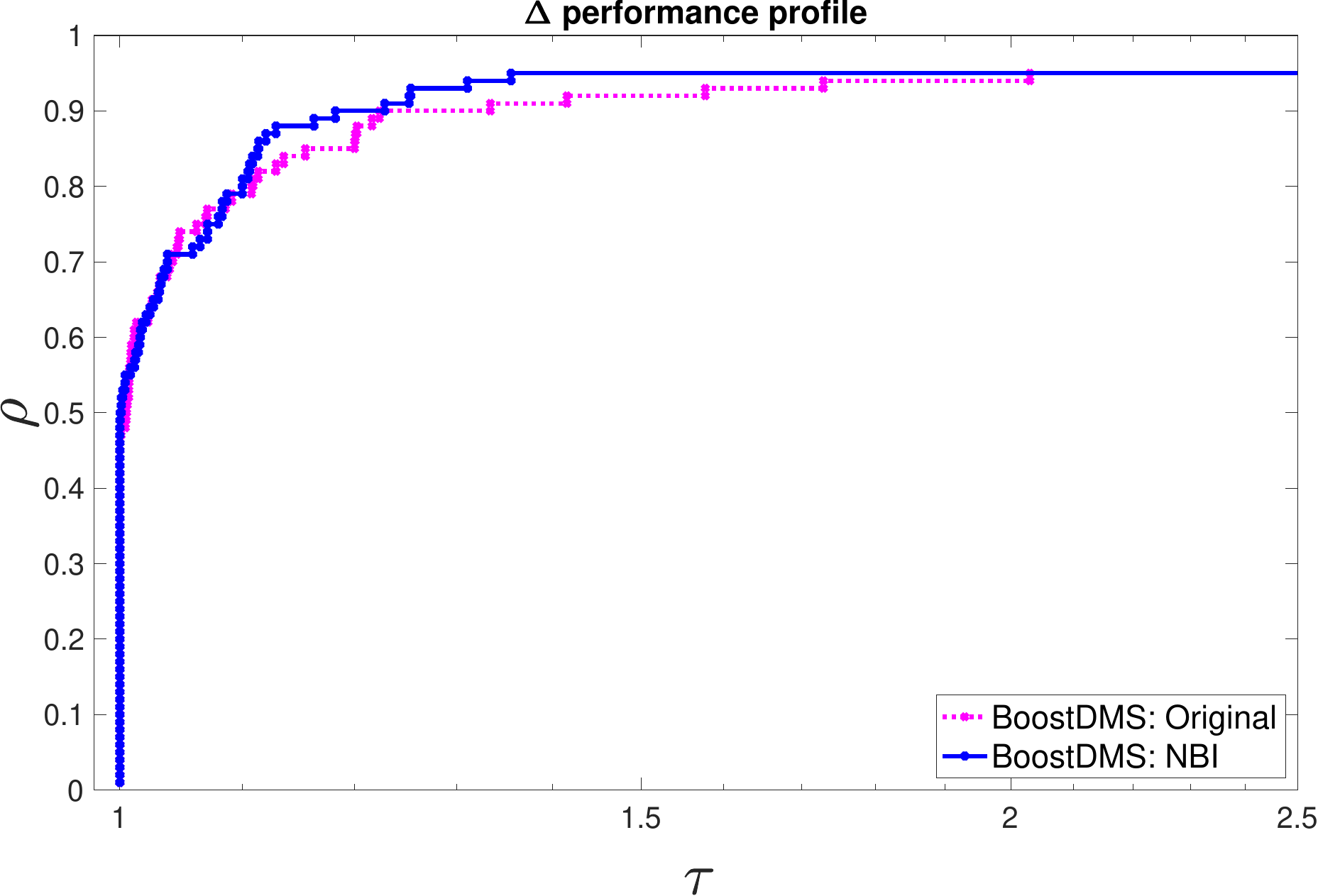}}
\caption{Comparing BoostDMS(Original) with BoostDMS(NBI) based on
performance profiles, for a maximum of $5000$ function
evaluations.} \label{BoostDMS_NBI}
\end{figure}

The NBI method is first compared with the original strategy
implemented in BoostDMS, based on the min-max scalarization.
Figure~\ref{BoostDMS_NBI} presents the performance profiles
obtained for the four metrics considered.

The profiles show that the NBI method behaves similarly to the original strategy in terms of purity and hypervolume. In both cases, the corresponding curves remain close, and no clear dominance of one strategy over the other can be observed. This indicates that incorporating NBI into the search step does not lead to a systematic improvement in the indicators more directly related to the quality of the final approximation, namely, the proportion of points belonging to the reference nondominated front and the dominated volume covered in the objective space.

A different behavior is observed for the spread metrics. For both
$\Gamma$ and $\Delta$ spreads, the curve associated with the NBI
method lies above the one corresponding to the original strategy.
This indicates that the use of NBI tends to generate
better-distributed approximations of the Pareto front. In
particular, the improvement observed in the $\Gamma$-spread
profile suggests a reduction in the largest gaps of the computed
front, while the $\Delta$-spread profile indicates a more uniform
distribution of nondominated points. This is consisting with the
original aims of NBI.

Overall, the main advantage of the NBI method appears to be
associated with the distribution of the final approximation rather
than with a clear improvement in purity or hypervolume. Within the
local model-based setting considered here, this advantage is
primarily reflected in the spread indicators, whereas the
quality-related metrics remain essentially comparable to those
obtained with the original strategy.

\subsection{Adapted $\epsilon-$constraint method}
\begin{figure}[H]
\centering
\subfigure[Purity]{\includegraphics[scale=0.25]{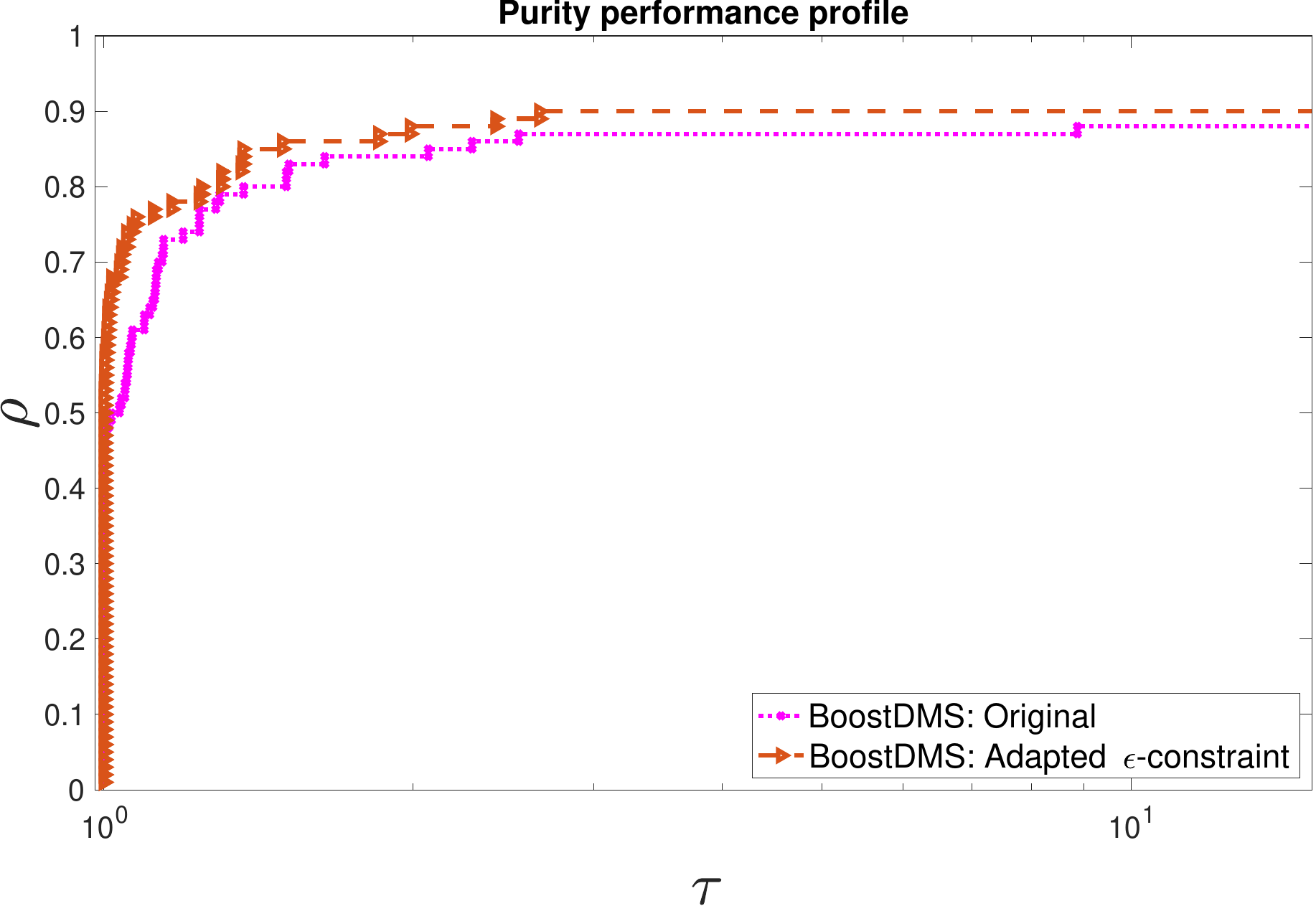}}
\subfigure[Hypervolume]{\includegraphics[scale=0.25]{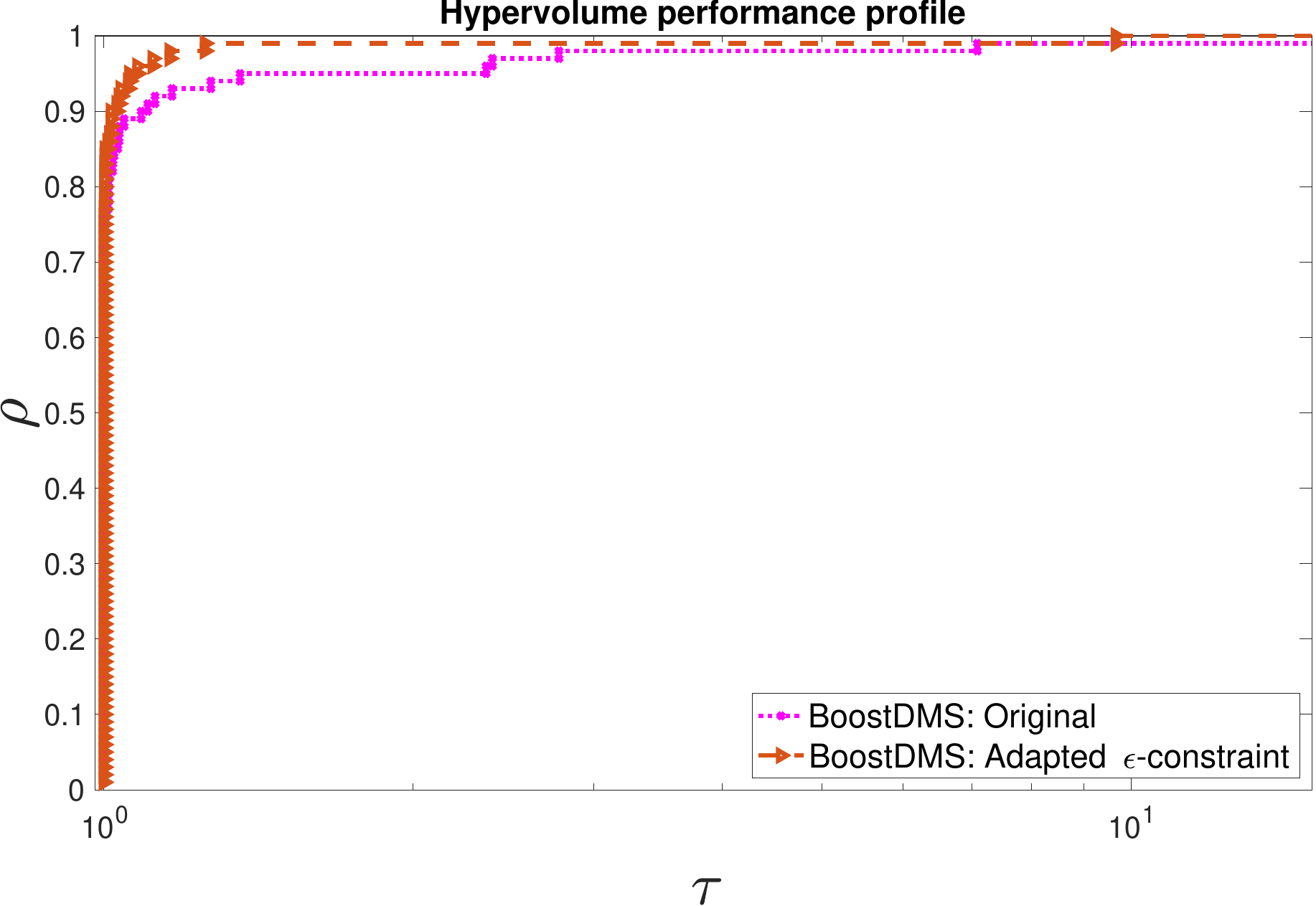}}
\subfigure[Spread Gamma
]{\includegraphics[scale=0.25]{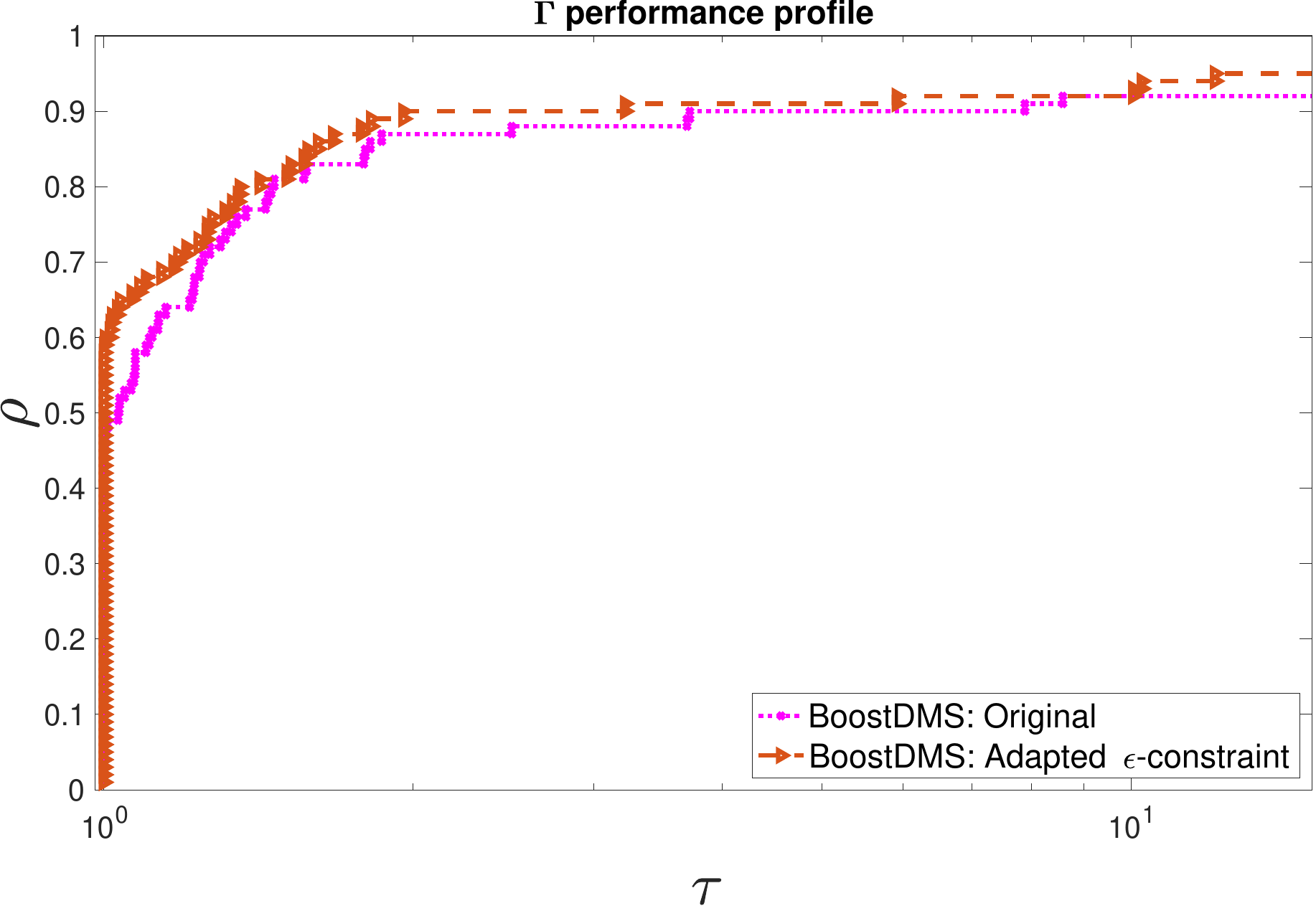}}
\subfigure[Spread
Delta]{\includegraphics[scale=0.25]{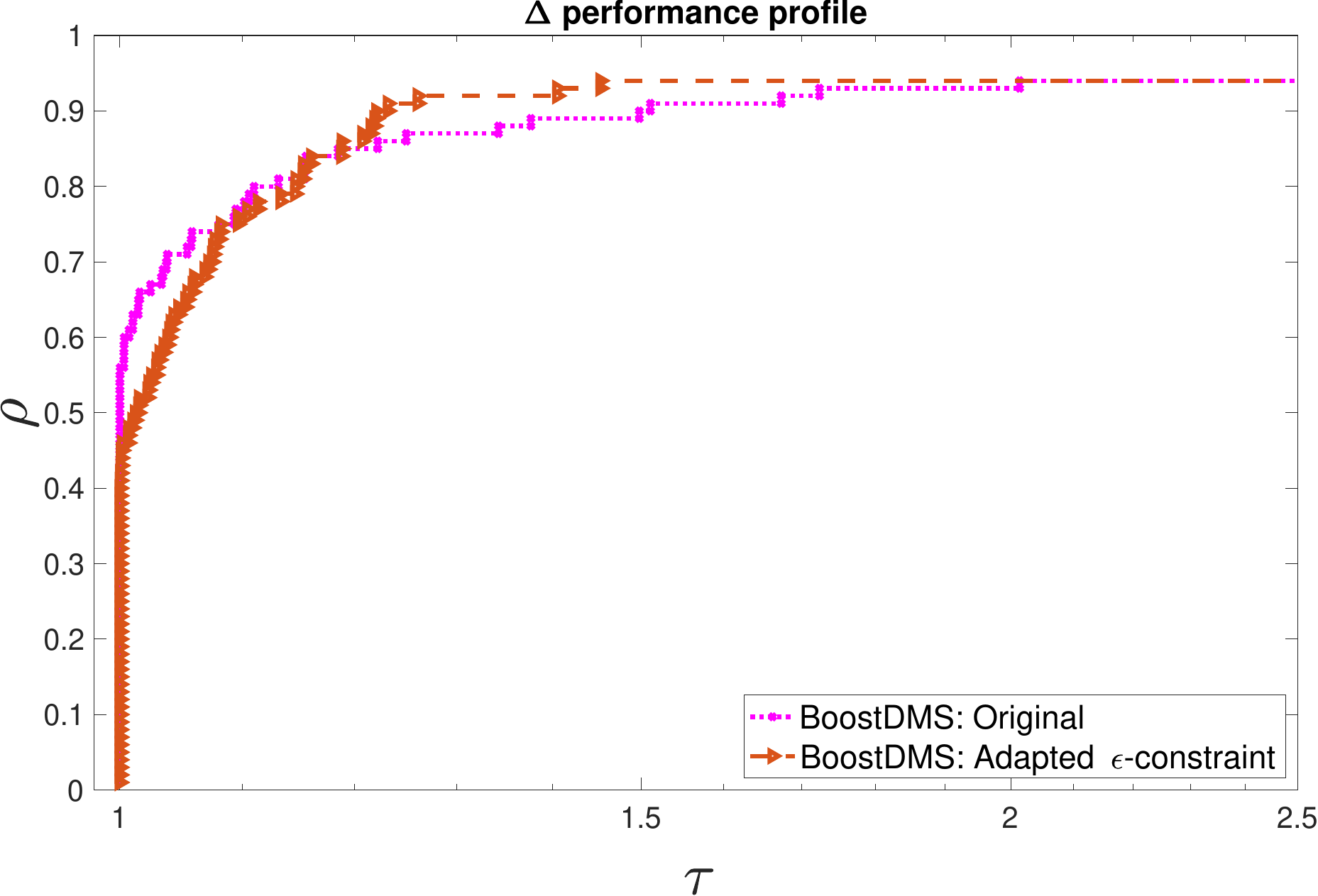}}
\caption{Comparing BoostDMS(Original) with BoostDMS(Adapted
$\epsilon-$constraint) based on performance profiles, for a
maximum of $5000$ function evaluations.}
\label{BoostDMS_EPSConstraint}
\end{figure}

The $\epsilon$-constraint method is now compared with the original
strategy implemented in BoostDMS.
Figure~\ref{BoostDMS_EPSConstraint} presents the performance
profiles obtained for the four metrics under analysis: purity,
hypervolume, $\Gamma$ spread, and $\Delta$ spread.

The most relevant differences are observed in the purity and hypervolume profiles. Regarding purity, the adapted $\epsilon$-constraint method outperforms the original strategy on a significant portion of the test set. This suggests that the final approximations generated by this method contain a larger proportion of points belonging to the reference nondominated front.

The hypervolume profile exhibits a similar trend. The curve
associated with the adapted $\epsilon$-constraint method is
generally above the one corresponding to the original strategy,
indicating that the computed approximations tend to dominate a
larger region of the objective space. The improvement is
noticeable and relatively consistent throughout the profile.

The comparison is more balanced for the spread indicators. For
$\Gamma$ spread, the adapted $\epsilon$-constraint method presents
a slight advantage, suggesting a reduction in the largest gaps
within the computed approximations. For $\Delta$ spread, the two
strategies exhibit similar behavior, with exception of $\tau=1$,
where there is an advantage of the original strategy. Even so, the
use of the adapted $\epsilon$-constraint method does not
substantially affect the uniformity of the final approximations,
remaining competitive with the original strategy in this regard.

Overall, the adapted $\epsilon$-constraint method primarily
improves the quality of the final approximations, as indicated by
the superior purity and hypervolume profiles.

\subsection{The IFSD method}
\begin{figure}[h!]
\centering
\subfigure[Purity]{\includegraphics[scale=0.25]{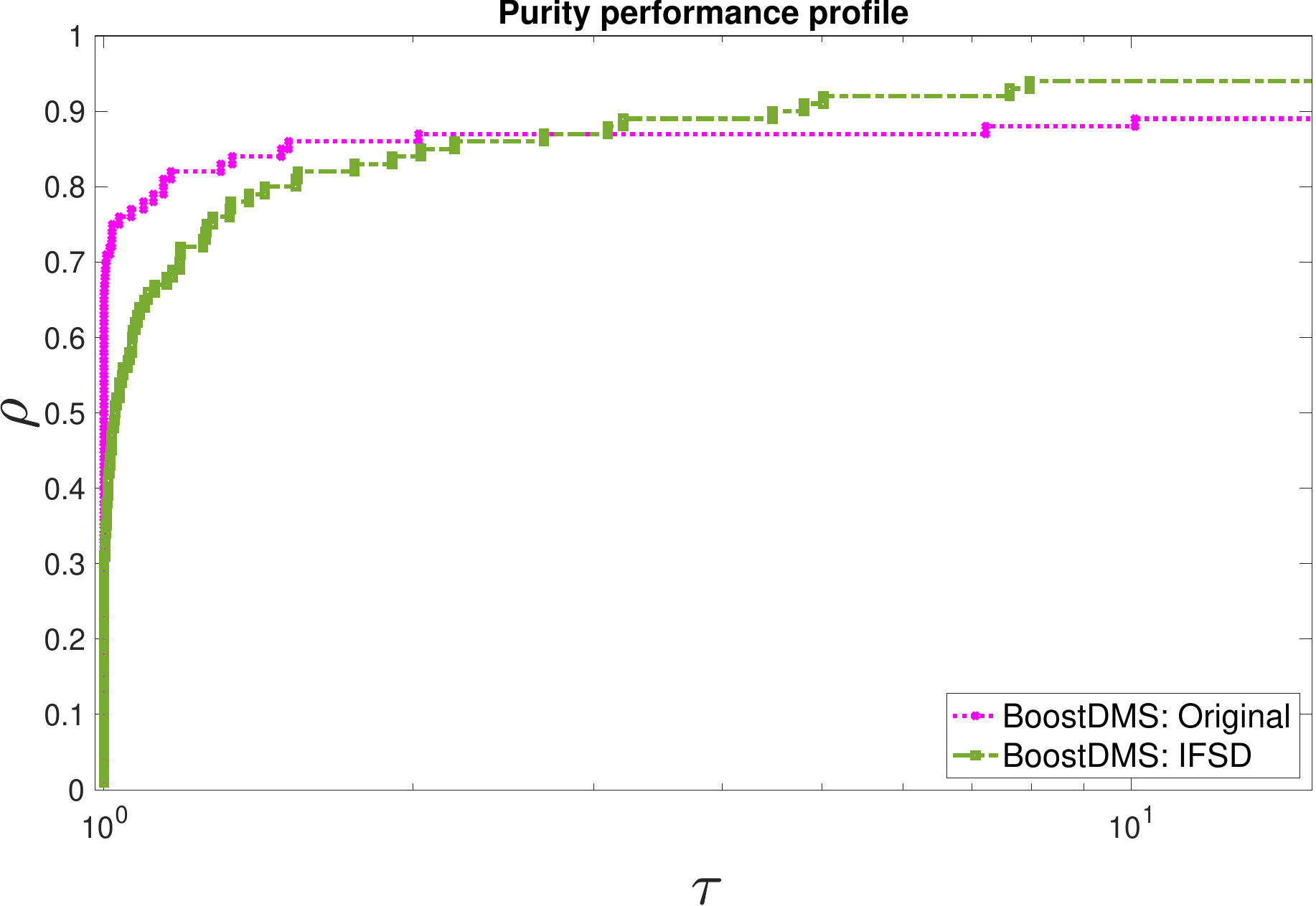}}
\subfigure[Hypervolume]{\includegraphics[scale=0.25]{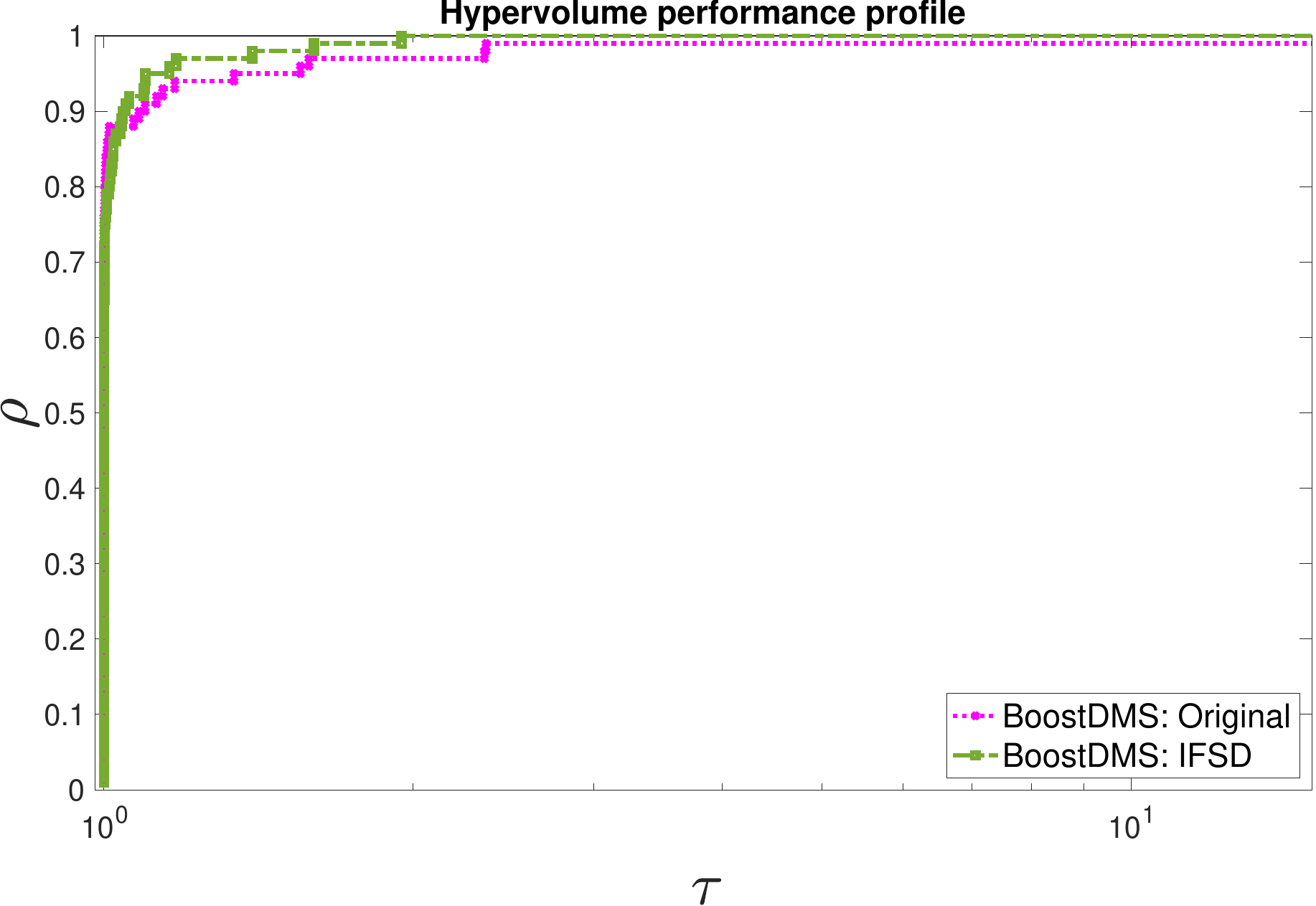}}
\subfigure[Spread Gamma
]{\includegraphics[scale=0.25]{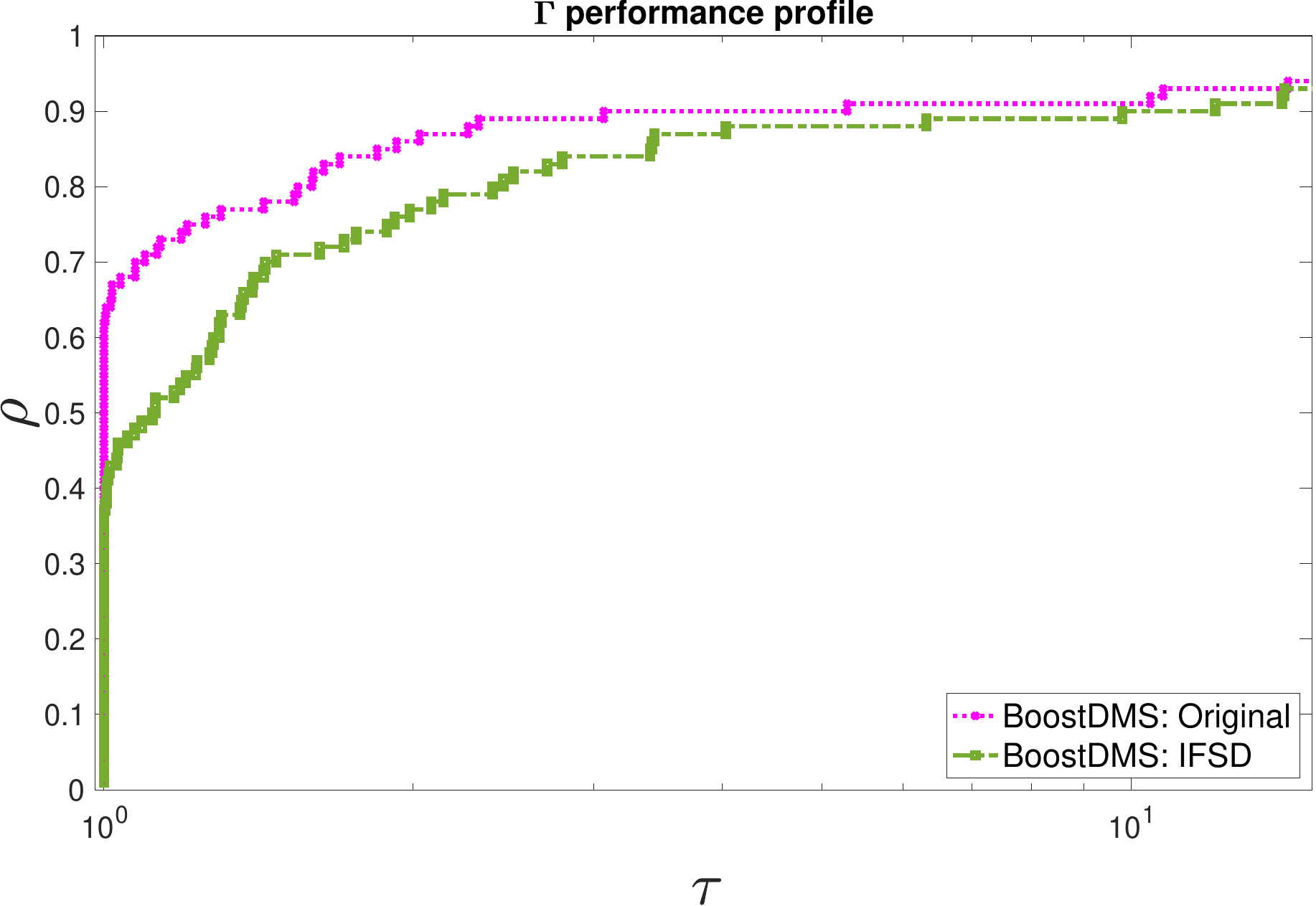}}
\subfigure[Spread
Delta]{\includegraphics[scale=0.25]{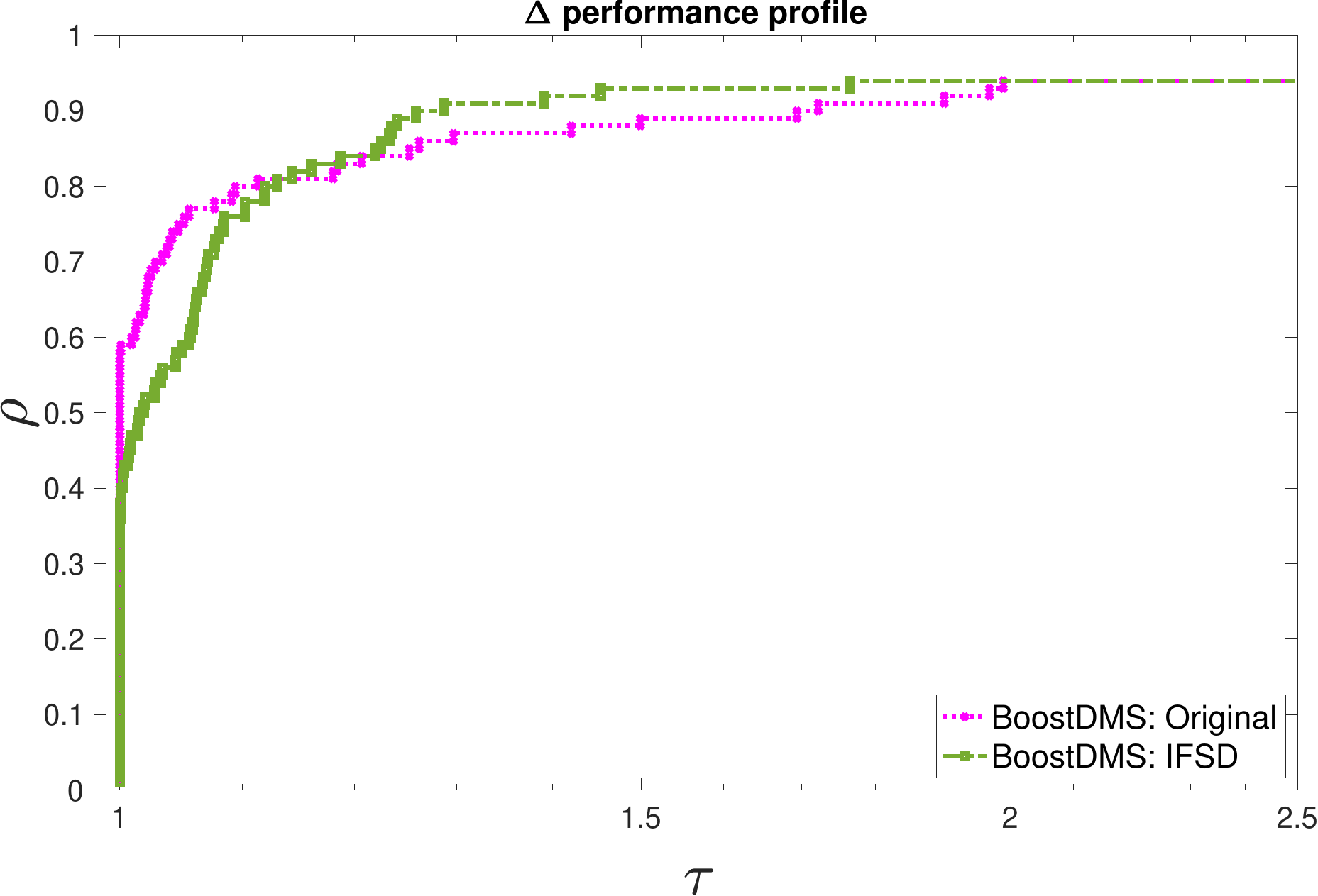}}
\caption{Comparing BoostDMS(Original) with BoostDMS(IFSD) based on
performance profiles, for a maximum of $5000$ function
evaluations.} \label{BoostDMS_IFSD}
\end{figure}

The final comparison is presented in Figure~\ref{BoostDMS_IFSD}, where the results obtained by BoostDMS(Original) and BoostDMS(IFSD) are reported.
The hypervolume performance profile shows that the two variants
behave very similarly, with a slight advantage of the IFSD
strategy. This indicates that, in terms of the overall quality of
the approximations (combining convergence and diversity), the IFSD
variant performs comparably to the original strategy. Moreover,
both variants exhibit high efficiency and robustness according to
this metric across the considered test set.

For purity, a different behavior is observed. BoostDMS(Original)
is the best performer for small $\tau$ values, indicating greater
efficiency. However, as $\tau$ increases, the curve associated
with BoostDMS(IFSD) eventually matches or surpasses that of
BoostDMS(Original), reaching an equal or higher $\rho$ value for
large $\tau$ values. This behavior suggests that BoostDMS(IFSD)
provides greater robustness in terms of the purity metric.

Regarding the spread indicators $\Delta$ and $\Gamma$,
BoostDMS(Original) is the most efficient variant, attaining the
minimum value for both metrics in $60\%$ of the problems.
BoostDMS(IFSD) is able to achieve a similar performance to the one
of BoostDMS(Original), but requiring larger values of $\tau$.

\section{Conclusions and future work}
\label{sec:conclusion} In this work, alternative strategies for
exploiting quadratic polynomial models within the search step of
Direct Multisearch were proposed and numerically assessed. The
considered approaches were based on the Normal Boundary
Intersection method, an adapted $\epsilon$-constraint
scalarization, and the Improved Front Steepest Descent algorithm.

The numerical results indicate that the different strategies
influence the quality of the final approximations in distinct
ways. The NBI-based search step produced results comparable to
those of the original BoostDMS strategy in terms of purity and
hypervolume, while providing some advantages in the spread
indicators. In contrast, the adapted $\epsilon$-constraint
strategy exhibited the most consistent improvement among the
tested alternatives. In particular, it achieved better performance
profiles for purity and hypervolume, while remaining competitive
with respect to the spread indicators.

The IFSD-based strategy exhibited a different behavior. Although it was designed to enrich the Pareto front approximation through refinement and exploration mechanisms, the numerical results did not reveal a clear advantage over the original strategy in terms of the quality-oriented indicators.

Overall, the results demonstrate that the choice of the
model-minimization strategy within the search step of DMS can
significantly influence the numerical performance of the
algorithm. Among the approaches evaluated, the adapted
$\epsilon$-constraint method appears to be the most promising,
particularly when the objective is to improve the quality of the
final approximation under a fixed budget of function evaluations.
These findings suggest that the model-minimization strategy should
be regarded as an algorithmic design choice rather than a
secondary implementation detail.

Several research directions may be explored to further enhance the numerical performance of BoostDMS. Promising extensions include the development of a search step based on the construction of quadratic models in low-dimensional subspaces at each iteration~\cite{Cartis_Roberts_2023}, as well as the incorporation of parallel computing techniques to improve the efficiency of the algorithm~\cite{Tavares_Bras_Custodio_Duarte_Medeiros_2023}.

{\small \printbibliography}


\end{document}